# Runge-Kutta symmetric interior penalty discontinuous Galerkin methods for modified Buckley-Leverett equations

Hong Zhang · Yunrui Guo · Weibin Li · Paul Andries Zegeling

**Abstract** We present a robust and accurate numerical method to solve the modified Buckley-Leverett equation in two-phase porous media flow with dynamic capillary pressure effect. A symmetric interior penalty discontinuous Galerkin method is used to discretize the equation in the space direction. For accuracy and stability issues, the third-order strong stability preserving implicit-explicit Runge-Kutta method is adopted to solve the nonlinear semi-discrete system: the linear diffusion term is discretized implicitly while the nonlinear flux term is discretized explicitly. The spatial accuracy of the discontinuous Galerkin method depends on the limiters applied to the solution: we test a minmod-TVB limiter, a simple WENO limiter and a high-order shock-capturing moment limiter to demonstrate that a suitable shock capturing moment limiter leads to more accurate approximation of solution. A set of representative numerical experiments are presented to show the accuracy and efficiency of the proposed approach. The results indicate that the moment limiter proposed by Moe et al. [Arxiv:1507.03024, 2015] is the most suitable one to be used in solving the modified Buckley-Leverett equation, and high order schemes perform much better than lower order schemes. Our simulation results are consistent with the previous results in Kao et al.[J. Sci. Comput., 64(3) (2015), 837-857], Zhang and Zegeling [J. Comput. Phys., 345 (2017), 510-527 and Commun. Comput. Phys., 22(4) (2017), 935-964].

## 1 Introduction

Two phase flow models in porous media have appeared in a wide range of areas, such as oil recovery, soil science and carbon sequestration etc. To illustrate the importance of such models, let us consider the example of water injecting into initially dry sandy porous media. Laboratory experiments [28,29,9] have revealed that non-monotone saturation profiles and finger patterns may be encountered under certain conditions. In the recent work [33,10], the authors discussed a modified Buckley-Leverett equation (MBLE) describing two-phase flow in porous media. This modification includes a third order mixed derivative term resulting from

Corresponding author: W. Li
H. Zhang
E-mail: h.zhang4@uu.nl
Department of Mathematics, Faculty of Science, Utrecht University, Budapestlaan 6, 3584CD Utrecht, The Netherlands

Y. Guo
E-mail: yunruiguo@nudt.edu.cn
State Key Laboratory of High Performance Computing, National University of Defense Technology, Changsha, Hunan 410073, China

W. Li
E-mail: liweibin@nudt.edu.cn
Computational Aerodynamics Institute, China Aerodynamics Research and Development Center, Mianyang 621000, China

P. A. Zegeling
E-mail: p.a.zegeling@uu.nl
Department of Mathematics, Faculty of Science, Utrecht University, Budapestlaan 6, 3584CD Utrecht, The Netherlands

the dynamic capillary pressure effect [14]. The dynamic capillary pressure coefficient serves as a bifurcation parameter which is critical in determining the type of solution profiles. In certain conditions, the solution may present non-monotone profile and even exhibit a damped oscillation.

To study the solution properties of the MBLE, a number of researches have been carried out recently. Most of them are in the finite difference and finite volume approaches. Peszynska and Yi [25] proposed a cell-centered finite difference method and a locally conservative Eulerian-Lagrangian method, but they noticed that such methods may cause instabilities in convection-dominated cases and for large dynamic effects. Van Duijn et al. [33] developed a finite difference method which combined a minmod slope limiter based on the first order upwind and Richtmyer's schemes. The obtained solutions agreed well with the traveling wave (TW) results. Wang and Kao [34] extended the second and third order central schemes to capture the nonclassical solutions of the MBLE. Kao further et al. [19] split the MBLE into a high-order linear equation and a nonlinear convective equation, and then integrated the linear equation with a pseudo-spectral method and the nonlinear equation with a Godunov-type central-upwind finite volume scheme. The computed solutions demonstrate that the higher-order spatial reconstruction using fifth-order WENO5 scheme gives more accurate numerical solutions. Hong et al. [16] adopted a fourth-order central difference scheme to resolve the spatial resolution and a standard fourth-order Runge-Kutta scheme to march the resulting algebraic system in time, they observed high wave number oscillatory waves under certain parametric conditions. Later work of de Moraes et al. [24] shows that those oscillatory waves do not satisfy threshold for the existence of non-monotonic wave fronts [33]. Thus they suggested to use schemes with nonlinear numerical stability properties to capture the different shock waves, as well as rarefaction waves. In [1], Abreu and Vieira discussed two numerical schemes based on the operator splitting technique. They found that the standard operator splitting may fail to capture the correct behavior of the solutions. In this sense, the operator splitting must take into account the dispersive-like character in both splitting steps. The authors also presented a non-splitting numerical method which is based on a fully coupled space-time mixed hybrid finite element/volume discretization approach. The numerical results suggest that such approach is able to account for the delicate nonlinear balance between the hyperbolic flux and the pseudo-parabolic term, but linked to a natural dispersive-like character of the full pseudo-parabolic differential equation.

Besides the standard finite difference and finite volume approaches, several adaptive mesh methods have also been proposed for the two-phase flow equations incorporating dynamic capillary pressure effect. Hu and Zegeling [17] applied a moving mesh finite element method to discretize the relaxation non-equilibrium Richards equation in the space direction. With the moving mesh technique, high mesh quality and accurate numerical solutions are obtained successfully. Refs. [35,37,36] studied the MBLE with adaptive moving mesh finite difference methods, their results show that to achieve the same accuracy, the adaptive methods need around a factor of 4-10 fewer grid points than the uniform grid cases.

The studies in [19] indicate that using high-order spatial reconstruction leads to more accurate approximation of solutions. To systematically study the effects of high-order schemes, we formulate and implement a discontinuous Galerkin (DG) discretization for the MBLE in two dimensional spaces. The objective is to investigate the applicability of DG method and the influences of different limiting strategies to the MBLE. Although several different DG methods [27,22] have been successfully applied to the standard Buckley-Leverett equation, there lacks research in the simulation of the MBLE. Moreover, it is far from clear that the method is necessarily able to capture shock profiles as well as rarefaction waves.

The remainder of this paper is organized in the following manner. Section 2 introduces the modified Buckley-Leverett equation, and the traveling wave results to be used in verifying the numerical solutions. Section 3 discusses the discretization of the equation by applying a symmetric interior penalty discontinuous Galerkin method in the spatial direction and an implicit-explicit Runge-Kutta method in the temporal direction. In Section 4, we present a minmod-TVB limiter, a WENO limiter and a shock-capturing moment limiter which will be used to limit the DG solutions. In section 5, we apply different initial conditions and parameters to the MBLE, to obtain shocks and rarefaction waves. Several 1D and 2D numerical experiments are carried out to demonstrate the effectiveness and advantages of the proposed scheme. Finally, Section 6 will outline the conclusions which are drawn from the results of the numerical schemes.

## 2 Mathematical model

In this section, we recall the two-phase flow equations and present some TW results. For more details see, e.g., [10].

2.1 The fractional flow formulation of two-phase flow equation with dynamic capillary pressure

Our model is for a homogeneous porous medium with a constant porosity $\phi$ and a constant intrinsic permeability $K$. Here we use the fractional flow formulation to describe two-phase wetting-non-wetting immiscible flow in two dimensional space. The saturation of each phase is defined as the volumetric fraction of the volume occupied by that phase. Denote the saturation of the wetting phase by $u$, then for a fully saturated porous medium, the saturation of the non-wetting phase is $1 - u$. The mass conservation equations for the two phases read

$$\frac{\partial(\phi\rho_w u)}{\partial t} + \nabla \cdot (\rho_w \mathbf{v}_w) = 0, \tag{1}$$

$$\frac{\partial(\phi\rho_n(1-u))}{\partial t} + \nabla \cdot (\rho_n \mathbf{v}_n) = 0, \tag{2}$$

in which $\rho_\alpha$ and $\mathbf{v}_\alpha, \alpha = n, w$ denote the density and the volumetric velocity of each phase.

Let the gravity acts in the negative $y$ axis, the Darcy's law reads

$$\begin{aligned}\mathbf{v}_\alpha &= -\frac{k_{r,\alpha}K}{\mu_\alpha}(\nabla p_\alpha + \rho_\alpha g \mathbf{e}_y), \\ &= -\lambda_\alpha(\nabla p_\alpha + \rho_\alpha g \mathbf{e}_y), \quad \alpha = n, w,\end{aligned} \tag{3}$$

where $g$ is the gravitational acceleration constant, $\mathbf{e}_y$ is the unit vector in the $y$ direction, $k_{r\alpha}$, $\mu_\alpha$, $p_\alpha$ and $\lambda_\alpha$ are the relative permeability function, viscosity, pressure and mobility of phase $\alpha$, respectively. Under non-equilibrium conditions, Stauffer [31], Hassanizadeh and Gray [14], Kalaydjian [18] proposed that the phases pressure difference $p_n - p_w$ can be written as a function of the equilibrium capillary pressure minus the product of the saturation rate of the wetting phase with a dynamic capillary coefficient $\tau$ [Pa s]:

$$p_n - p_w = P_c(u) - \tau \frac{\partial u}{\partial t}, \tag{4}$$

where $P_c$ modeling the capillary pressure - saturation relationship under an equilibrium condition, is a smooth and decreasing function of saturation $u$, and $\tau$ can be explained as a relaxation time. We refer to [13] for a review of experimental work on dynamic effects in the pressure-saturation relationship.

When the two phases (e.g. water and oil) are incompressible, by defining the total velocity $\mathbf{v}_T = \mathbf{v}_n + \mathbf{v}_w = [v_T^x, v_T^y]^T$ and the fractional flow rate of the wetting phase $f_w(u) = \frac{\lambda_w}{\lambda_w + \lambda_n}$, the velocity of the wetting phase can then be expressed by

$$\mathbf{v}_w = f_w(u)[v_T + \lambda_n(\nabla(p_n - p_w) - (\rho_w - \rho_n)g)]. \tag{5}$$

Substituting (5) into (1) and incorporating (4), we get a modified two-phase Buckley-Leverett equation,

$$\frac{\partial u}{\partial t} + \frac{\partial}{\partial x}F(u) + \frac{\partial}{\partial y}G(u) + \nabla \cdot [D(u)\nabla u] - \tau \nabla \cdot [H(u)\nabla \frac{\partial u}{\partial t}] = 0, \tag{6}$$

where

$$\begin{aligned} F(u) &= \frac{1}{\phi}f_w(u)v_T^x, & G(u) &= \frac{1}{\phi}f_w(u)[v_T^y - \lambda_n(u)(\rho_w - \rho_n)g], \\ D(u) &= \frac{1}{\phi}\lambda_n(u)f_w(u)P_c'(u), & H(u) &= \frac{1}{\phi}\lambda_n(u)f_w(u). \end{aligned} \tag{7}$$

In realistic modeling, the parameters and functions in (6) and (7) depend on the properties of the porous medium and the phases. In this work, we are interested in the relationship between the saturation $u$ and the dynamic coefficient $\tau$, thus we will consider a simplified MBLE with initial and boundary conditions

$$\begin{cases} \frac{\partial u}{\partial t} + \nabla \cdot \mathbf{v}(u) - \epsilon \Delta u - \tau \epsilon^2 \Delta(\frac{\partial u}{\partial t}) = 0, \\ u(x, y, t = 0) = u_0(x, y), \quad (x, y) \in \Omega \\ u(x, y, t) = u^D(x, y, t), \quad (x, y) \in \Gamma^D \\ \mathbf{n} \cdot \nabla u = 0, \quad (x, y) \in \Gamma^N, \end{cases} \tag{8}$$

where $\mathbf{v} = [F(u), G(u)]^T$ is the flux vector, $\Omega \subset R^2$ is the physical domain and $\Gamma^D$ is the Dirichlet boundary, $\Gamma^N$ is the Neumann boundary, $\mathbf{n}$ denotes the unit outward normal of $\partial \Omega$.

**Table 1** TW results of the 1D MBLE (9) summarized from Ref. [33].

| Region | Solution description |
|---|---|
| $(u_B, \tau) \in A_1$ | Rarefaction wave from $u_B$ down to $u_\alpha$ trailing an admissible Lax shock from $u_\alpha$ down to $u_0$ |
| $(u_B, \tau) \in A_2$ | Rarefaction wave from $u_B$ down to $\bar{u}$ trailing an undercompressive shock from $\bar{u}$ down to $u_0$ |
| $(u_B, \tau) \in B$ | An admissible Lax shock from $u_B$ up to $\bar{u}$ (may exhibit oscillations near $u_+ = u_B$) trailing an undercompressive shock from $\bar{u}$ down to $u_0$ |
| $(u_B, \tau) \in C_1$ | An admissible Lax shock from $u_B$ down to $u_0$ |
| $(u_B, \tau) \in C_2$ | An admissible Lax shock from $u_B$ down to $u_0$ (may exhibit oscillations near $u_+ = u_B$) |

2.2 Traveling wave solutions for one dimensional MBLE

The 1D MBLE in the $x$-direction reads

$$\frac{\partial u}{\partial t} + \frac{\partial F(u)}{\partial x} - \epsilon \frac{\partial^2 u}{\partial x^2} - \tau \epsilon^2 \frac{\partial^3 u}{\partial x^2 \partial t} = 0. \tag{9}$$

Consider a TW solution connecting $u_-$ and $u_+$ ($u_+ > u_- \in [0,1]$), by introducing the TW coordinate $\eta = x - st$ and substituting the TW solution $u(\eta)$ into (9) we obtain a third order ordinary differential equation (ODE)

$$\begin{cases} -su' + [F(u)]' - \epsilon u'' + s\tau\epsilon^2 u''' = 0, \\ u(\pm\infty) = u(\pm), \quad u_+ > u_- \in [0,1], \\ u'(\pm\infty) = u''(\pm\infty) = 0, \end{cases}$$

where prime denotes differentiation with respect to $\eta$, the boundary conditions of the ODE are obtained by the definition of TW solutions. Integrating this equation over $(-\infty, \eta)$ yields the second-order ODE:

$$\begin{cases} -s(u - u_-) + [F(u) - F(u_-)] - \epsilon u' + s\tau\epsilon^2 u'' = 0, \\ u(\pm\infty) = u_\pm, \end{cases} \tag{10}$$

with $s$ determined by the Rankine-Hugoniot condition

$$s = \frac{F(u_+) - F(u_-)}{u_+ - u_-}.$$

Consider $F(u) = \frac{u^2}{u^2 + 0.5(1-u)^2}$, $\epsilon = 10^{-3}$, for a fixed value of $u_-$, the dependency between $\tau$ and the value $u_+$ has been analyzed in Ref. [33, 30, 32]. Let $u_I$ be the unique inflection point of the flux function $F(u)$, we summarize the results obtained by [33].

When $u_0 \in [0, u_I)$, it is proved that there is a constant $\tau_*$ such that for all $\tau \in [0, \tau_*]$, there exists a unique solution of (10) connecting $u_+ = u_\alpha$ and $u_- = u_0$, where $u_\alpha$ is the unique root of the equation

$$F'(u) = \frac{F(u) - F(u_0)}{u - u_0}.$$

When $\tau > \tau_*$, there exists a unique constant $\bar{u} > u_\alpha$, such that (10) has a unique solution connecting $u_+ = \bar{u}$ and $u_- = u_0$. For $u_- = u_0 < u_+ = u_B < \bar{u}(\tau)$, the solution of (10) will exist only if $u_B \in (u_0, \underline{u})$, where $\underline{u}$ is the unique root in the interval $(u_0, \bar{u})$ of

$$\frac{F(u) - F(u_0)}{u - u_0} = \frac{F(\bar{u}) - F(u_0)}{\bar{u} - u_0}.$$

When $\tau > \tau_*$ and $u_B \in (\underline{u}, \bar{u})$, there is no TW solution of (10) connecting $u_+ = u_B$ and $u_- = u_0$. In this situation, the solution profile is non-monotonic, two TWs are used in succession: one from $u_+ = u_B$ to $u_- = \bar{u}$ and one from $u_+ = \bar{u}$ to $u_- = u_0$. For any $u_B \in (\underline{u}, \bar{u})$ and $\tau > \tau_*$, there exists a unique solution of (10) such that $u_+ = u_B$, $u_- = \bar{u}$.

When $u_0 < u_I$ and $u_B > u_0$, the traveling solutions can be classified using the five regions in the bifurcation diagram. The results summarized from Ref. [33] are given in Table 1, and the bifurcation diagram for $u_0 = 0$ is plotted in Fig. 1.

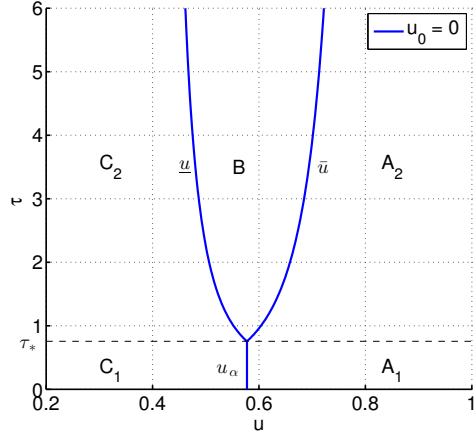

**Fig. 1** Bifurcation diagram with $u_0 = 0$.

## 3 The discontinuous Galerkin method with implicit-explicit time-marching

In this section, we describe the discontinuous Galerkin discretization and the implicit-explicit (IMEX) Runge-Kutta scheme for the governing equation (8) in detail.

3.1 Semi-discrete formulation

In the DG method, the 2D physical domain $\Omega$ is subdivided into a set of rectangular elements $\Omega_h = \{E_{i,j}\}$. Let $\Gamma_h$, $\Gamma_h^D$ and $\Gamma_h^N$ be the unions of the interior faces in $\Omega_h$, the Dirichlet and Neumann boundary interfaces, respectively. For any interface $e = \partial E^- \cap \partial E^+$, we let $\mathbf{n}_e$ pointing from $E^-$ to $E^+$ denote the unit outward normal to the face $e$ of $E^-$. We also denote by $|e|$ the measure of $e$. The DG finite element space $V_h^k$ is given by

$$V_h^k = \{v(x,y) \in L^2(\Omega) : v|_E \in P^k(E), \forall E \in \Omega_h\}, \tag{11}$$

where $P^k(E)$ is the set of polynomials defined on element $E$ of degree at most $k$. Notice that the piece-wise polynomial $v \in V_h^k$ is allowed to be discontinuous at the cell interface $e$. Thus, we define the jump operator $[\cdot]$ and the average operator $\{\cdot\}$ along the interface $e$ as

$$\forall e \in \Gamma_h, \qquad [v] := v^- - v^+, \quad \{v\} := \frac{1}{2}(v^- + v^+),$$
$$\forall e \in \Gamma_h^D \cup \Gamma_h^N, \quad [v] := v^-, \qquad \{v\} := v^-,$$

where $v^- := v|_{E^-}, v^+ := v|_{E^+}$ are the traces of $v$ on elements $E^-$ and $E^+$.

Let $\{\psi_l\}_{l=0,\cdots,N_k-1}$ denote a local polynomial basis of $V_h^k$ on element $E$, an approximations to the solution $u$ has the form

$$u_h(x,y,t)|_E = \sum_{l=0}^{N_k-1} u_j(t)\psi_l(x,y), \quad \forall (x,y) \in E. \tag{12}$$

The variational formulation for the MBLE can be derived using the Galerkin approach. First, multiply (8) by a smooth test function $\varphi : \Omega \to R$ and then integrate by parts over each element. Hence, the symmetric interior penalty discontinuous Galerkin (SIPDG) discretization to (8) is described as: find $u_h \in V_h^k$ such that $\forall \varphi \in V_h^k$,

$$<\frac{\partial u_h}{\partial t}, \varphi>_E + \tau\epsilon^2 \mathcal{A}_{diff,1}(\frac{\partial u_h}{\partial t}, \varphi) + \epsilon \mathcal{A}_{diff,2}(u_h, \varphi) + \mathcal{A}_{adv}(u_h, \varphi) = r(\varphi), \quad \forall \varphi \in V_h^k, \forall E \in \Omega_h, \tag{13}$$

where $< \cdot, \cdot >_E$ denotes the $L^2(\Omega_h)$ inner product on element $E$. The bilinear operators $\mathcal{A}_{diff,1}$, $\mathcal{A}_{diff,2}$ and $\mathcal{A}_{adv}$ are expressed as

$$\mathcal{A}_{diff,1}(\frac{\partial u_h}{\partial t}, \varphi) = (\sum_{E \in \Omega_h} \int_E \nabla \frac{\partial u_h}{\partial t} \cdot \nabla \varphi \mathrm{d}x - \sum_{e \in \Gamma_h} \int_e \{\nabla \frac{\partial u_h}{\partial t} \cdot \mathbf{n}_e\}[\varphi]\mathrm{d}s - \sum_{e \in \Gamma_h} \int_e \{\nabla \varphi \cdot \mathbf{n}_e\}[\frac{\partial u_h}{\partial t}]\mathrm{d}s$$
$$+ \sum_{e \in \Gamma_h} \sigma \int_e [\frac{\partial u_h}{\partial t}][\varphi]\mathrm{d}s), \qquad (14)$$

$$\mathcal{A}_{diff,2}(u_h, \varphi) = (\sum_{E \in \Omega_h} \int_E \nabla u_h \cdot \nabla \varphi \mathrm{d}s - \sum_{e \in \Gamma_h \cup \Gamma_h^D} \int_e \{\nabla u_h \cdot \mathbf{n}_e\}[\varphi]\mathrm{d}s - \sum_{e \in \Gamma_h \cup \Gamma_h^D} \int_e \{\nabla \varphi \cdot \mathbf{n}_e\}[u_h^n]\mathrm{d}s$$
$$+ \sum_{e \in \Gamma_h \cup \Gamma_h^D} \sigma \int_e [u_h][\varphi]\mathrm{d}s), \qquad (15)$$

$$\mathcal{A}_{adv}(u_h, \varphi) = -\sum_{E \in \Omega_h} \int_E \mathbf{v}(u_h) \cdot \nabla \varphi \mathrm{d}x + \sum_{e \in \Gamma_h} \int_e \mathbf{v}(u_h) \cdot \mathbf{n}_e \varphi \mathrm{d}s. \qquad (16)$$

In (16) the flux in the direction of the out unit normal $\mathbf{v}(u_h) \cdot \mathbf{n}$ is not defined on the interface $e$ because of the discontinuity of $u_h$ across the interface. Therefore, we approximate the normal trace by the local Lax-Friedrichs flux function

$$\mathbf{v}^*(u^-, u^+) = \frac{1}{2}\mathbf{n}_e \cdot [\mathbf{v}(u^-) + \mathbf{v}(u^+)] - \frac{c}{2}\mathbf{n}_e(u^+ - u^-), \quad c = \max_{u \in [\min(u^-, u^+), \max(u^-, u^+)]} |\mathbf{n}_e \cdot \frac{\partial \mathbf{v}}{\partial u}|. \qquad (17)$$

Using the numerical flux (17), Eq. (16) is then replaced by

$$\mathcal{A}_{adv}(u_h, \varphi) = -\sum_{E \in \Omega_h} \int_E \mathbf{v}(u_h) \cdot \nabla \varphi \mathrm{d}x + \sum_{e \in \Gamma_h} \int_e \mathbf{v}^*(u^-, u^+) \cdot \mathbf{n}_e \varphi \mathrm{d}s, \qquad (18)$$

The boundary conditions are incorporated in the face integral in the right hand side term $r(\varphi)$:

$$r(\varphi) = -\sum_{e \in \Gamma_h^D} \int_e \mathbf{v}(u^D) \cdot \mathbf{n} \varphi \mathrm{d}s - \epsilon(\sum_{e \in \Gamma_h^D} \int_e u^D \{\nabla \varphi \cdot \mathbf{n}_e\}\mathrm{d}s + \sum_{e \in \Gamma_h^D} \sigma \int_e u^D [\varphi]\mathrm{d}s). \qquad (19)$$

The parameter $\sigma > 0$ in (14) and (15) denotes the penalty, which is chosen as

$$\sigma_{e,E} = k(k+1)\frac{|e|}{|E|}. \qquad (20)$$

If the face $e$ is at the boundary, choose $\sigma = \sigma_{e,E}$. For an interior face, we take the average of the two values at this face. Similar choice of $\sigma$ has also been used in [4].

Testing (13) with $\varphi = \psi_l$ for $l \in [0, \cdots, N_k - 1]$ and substituting (12) into (13) will yield a time-dependent system of equations. Written in matrix from, the system is then given by

$$(M + \tau\epsilon^2 A_{diff,1})\bar{u}_t + \epsilon A_{diff,2}\bar{u} + A_{adv}(u_h) = \bar{r}(\varphi), \forall E \in \Omega_h, \qquad (21)$$

where $\bar{u}$ is a representation vector of the coefficients $\{u_l\}_{l=0,\cdots,N_k-1}$ in (12).

3.2 Legendre basis

We implement the DG method using the unit Legendre basis defined on a reference element. In 2D, the physical element $E_{ij} = [x_{i-\frac{1}{2}}, x_{i+\frac{1}{2}}] \times [y_{j-\frac{1}{2}}, y_{j+\frac{1}{2}}]$ is mapped to the reference element $E_r = [0,1]^2$ via coordinates mapping

$$\xi = \frac{x - x_i}{\Delta x_i} + \frac{1}{2}, \quad \eta = \frac{y - y_j}{\Delta y_j} + \frac{1}{2}, \quad (x, y) \in E_{ij}. \qquad (22)$$

At any moment $t$, the DG solution $u_h(x, y, t)$ on the element $E_{ij}$ can be expressed in the basis function space $\{\psi_l^r(\xi, \eta)\}_{l=0,\cdots N_k-1}$ defined on $E_r$,

$$u_h(x, y, t)|_{E_{ij}} = u_h(\xi, \eta, t)|_{E_r} = \sum_{l=0}^{N_k-1} u_l(t)\psi_l^r(\xi(x), \eta(y)), \qquad (23)$$

where the degrees of freedom $\{u_l\}_{l=0,\cdots,N_k-1}$ are the so called modal coefficients. We point out that for this set of basis functions the average value of the solution on element $E_{ij}$ is $u_0$.

In the following we present the 1D and 2D basis functions defined on the reference element $[0,1]^{\text{DIM}}$.

1. In one dimensional space (DIM = 1), the unit Legendre basis functions defined on 1D interval $[0,1]$ read:

– $k = 1$:

$$\psi_0(\xi) = 1, \quad \psi_1(\xi) = 2\sqrt{3}(\xi - \frac{1}{2}),$$

– $k = 2$:

$$\psi_0(\xi) = 1, \quad \psi_1(\xi) = 2\sqrt{3}(\xi - \frac{1}{2}), \quad \psi_2(\xi) = \frac{\sqrt{5}}{2}(12(\xi - \frac{1}{2})^2 - 1),$$

– $k = 3$:

$$\psi_0(\xi) = 1, \quad \psi_1(\xi) = 2\sqrt{3}(\xi - \frac{1}{2}), \quad \psi_2(\xi) = \frac{\sqrt{5}}{2}(12(\xi - \frac{1}{2})^2 - 1), \psi_3(\xi) = \frac{\sqrt{7}}{2}(40(\xi - \frac{1}{2})^3 - 6(\xi - \frac{1}{2})).$$

2. In two dimensional space (DIM = 2), the unit Legendre basis functions defined on 2D square $[0,1]^2$ read:

– $k = 1$:

$$\psi_0(\xi,\eta) = 1, \quad \psi_1(\xi,\eta) = 2\sqrt{3}(\xi - \frac{1}{2}), \quad \psi_2(\xi,\eta) = 2\sqrt{3}(\eta - \frac{1}{2}),$$

– $k = 2$:

$$\psi_0(\xi,\eta) = 1, \quad \psi_1(\xi,\eta) = 2\sqrt{3}(\xi - \frac{1}{2}), \quad \psi_2(\xi,\eta) = \frac{\sqrt{5}}{2}(12(\xi - \frac{1}{2})^2 - 1),$$
$$\psi_3(\xi,\eta) = 2\sqrt{3}(\eta - \frac{1}{2}), \quad \psi_4(\xi,\eta) = 12(\xi - \frac{1}{2})(\eta - \frac{1}{2}), \quad \psi_5(\xi,\eta) = \frac{\sqrt{5}}{2}(12(\eta - \frac{1}{2})^2 - 1).$$

– $k = 3$:

$$\psi_0(\xi,\eta) = 1, \quad \psi_1(\xi,\eta) = 2\sqrt{3}(\xi - \frac{1}{2}), \quad \psi_2(\xi,\eta) = \frac{\sqrt{5}}{2}(12(\xi - \frac{1}{2})^2 - 1),$$
$$\psi_3(\xi,\eta) = \frac{\sqrt{7}}{2}(40(\xi - \frac{1}{2})^3 - 6(\xi - \frac{1}{2})), \quad \psi_4(\xi,\eta) = 2\sqrt{3}(\eta - \frac{1}{2}), \quad \psi_5(\xi,\eta) = 12(\xi - \frac{1}{2})(\eta - \frac{1}{2}),$$
$$\psi_6(\xi,\eta) = \sqrt{15}(12(\xi - \frac{1}{2})^2 - 1)(\eta - \frac{1}{2}), \quad \psi_7(\xi,\eta) = \frac{\sqrt{5}}{2}(12(\eta - \frac{1}{2})^2 - 1),$$
$$\psi_8(\xi,\eta) = \sqrt{15}(12(\eta - \frac{1}{2})^2 - 1)(\xi - \frac{1}{2}), \quad \psi_9(\xi,\eta) = \frac{\sqrt{7}}{2}(40(\eta - \frac{1}{2})^3 - 6(\eta - \frac{1}{2})).$$

3.3 Implicit-explicit Runge-Kutta time integration

The SIPDG discretization results in a large system of ODEs containing both stiff and nonstiff parts, which is suitable to be integrated using the implicit-explicit Runge-Kutta (IMEX RK) method. Here we illustrate the IMEX RK scheme by applying it to an ODE.

Consider the ODE

$$u_t = f(u) + g(u), \tag{24}$$

| | 0 | 0 | 0 | | 0 | 0 | 0 |
|---|---|---|---|---|---|---|---|
| 0 | | | | 0 | | | |
| 1 | 1 | 0 | 0 | 1 | $\frac{14}{15}$ | $\frac{1}{15}$ | 0 |
| $\frac{1}{2}$ | $\frac{1}{4}$ | $\frac{1}{4}$ | 0 | $\frac{1}{2}$ | $\frac{7}{30}$ | $\frac{1}{5}$ | $\frac{1}{15}$ |
| A | $\frac{1}{6}$ | $\frac{1}{6}$ | $\frac{2}{3}$ | $\tilde{A}$ | $\frac{1}{6}$ | $\frac{1}{6}$ | $\frac{2}{3}$ |

**Table 2** Butcher tableau for SSP3(3, 3, 3) method.

where $f(u)$ is a nonlinear term and $g(u)$ is a stiff term. The IMEX RK scheme treats the stiff term implicitly and the nonstiff term explicitly, thus reducing the computational complexity of the scheme because we do not need to limit the time step to satisfy the stability restriction for the stiff term.

We will adopt a strong stability preserving (SSP) IMEX RK scheme. An SSP IMEX RK scheme is referred to as SSP$k(s, \sigma, p)$, where $k$ is the order of the method in the stiff limit $\epsilon \to 0$, $s$ is the number of stages in the implicit scheme and $\sigma$ is the number of stages in the explicit scheme, $p$ is the global order of the resulting combined method. The SSP3(3, 3, 3) scheme proposed by [15] is presented in Table 2.

In the IMEX scheme, because of the explicit treating of the convection term, we let the time step $\Delta t$ satisfy the CFL condition

$$\Delta t = \text{CFL} \cdot \min\bigl(\frac{\Delta x}{|F'(u)|}, \frac{\Delta y}{|G'(u)|}\bigr), \tag{25}$$

where CFL is the Courant-Friedrichs-Levy constant, $F'(u)$ and $G'(u)$ are the wave speeds in $x$ and $y$ directions, respectively. In practice, as suggested by Cockburn et al. [8], for DG discretizations using polynomials of degree $k$, we take the CFL number as

$$\text{CFL} = \frac{1}{2k+1}. \tag{26}$$

## 4 Limiting strategies

The studies in Refs. [36,37] show that the MBLE with certain initial conditions and given values of $\tau$ will develop overshoots and sharp gradients along the moving front, which may cause significant difficulties in numerical simulations. In order to mitigate numerical oscillations and to preserve the original accuracy, we will apply limiters after each complete Runge-Kutta time step.

Many different types of limiters have been developed for DG method in the past several decades. For example, the minmod type total variation bounded (TVB) limiters [7], the moment based limiters [5] and the WENO limiters [26,40,11]. Since our objective is to accurately compute the overshoot values, we will try different limiters and find the one which can not only preserve the high order accuracy in smooth regions but also reduce the spurious oscillations near shocks or discontinuities.

In the following subsections, we present details of three limiting procedures for the DG method. To apply these limiters, we adopt the framework proposed by [26]:

1. Identify troubled cells which might need limiting procedure;
2. Replace the solution polynomials in those troubled cells with the reconstructed polynomials which maintain the original cell averages.

4.1 The minmod-TVB limiter

In one dimensional space, the minmod-TVB limiter limits the first order component of the solution on $E_i$ by using cell averages of the two neighboring cells $E_{i-1}$ and $E_{i+1}$. We denote the cell average of the solution $u_i(x)$ on $E_i$ as

$$\bar{u}_i = \frac{1}{\Delta x_i} \int_{E_i} u_i(x) \mathrm{d}x, \tag{27}$$

and the forward and backward differences as

$$\tilde{u}_i = u^-_{i+\frac{1}{2}} - \bar{u}_i, \quad \tilde{\tilde{u}}_i = \bar{u}_i - u^+_{i-\frac{1}{2}},$$
$$\Delta_- \bar{u}_i = \bar{u}_i - \bar{u}_{i-1}, \quad \Delta_+ \bar{u}_i = \bar{u}_{i+1} - \bar{u}_i.$$

We modified the slopes using the minmod-TVB limiter $\tilde{m}()$:

$$\tilde{u}_i^m = \tilde{m}(\tilde{u}_i, \frac{1}{2}\Delta_+\bar{u}_i, \frac{1}{2}\Delta_-\bar{u}_i), \quad \tilde{\tilde{u}}_i^m = \tilde{m}(\tilde{\tilde{u}}_i, \frac{1}{2}\Delta_+\bar{u}_i, \frac{1}{2}\Delta_-\bar{u}_i),$$
$$\tilde{m}(a_1, a_2, \cdots, a_n) = \begin{cases} a_1, & \text{if } |a_1| \leq M_{TVB}h^2, \\ m(a_1, a_2, \cdots, a_n), & \text{otherwise,} \end{cases} \quad (28)$$

where the TVB parameter $M_{TVB}$ has to be chosen adequately [7], and the minmod function $m()$ is

$$m(a_1, a_2, \cdots, a_n) := \begin{cases} \min(a_1, a_2, \cdots, a_n), & \text{if } a_i > 0, \forall i, \\ \max(a_1, a_2, \cdots, a_n), & \text{if } a_i < 0, \forall i, \\ 0, & \text{otherwise.} \end{cases} \quad (29)$$

The trace values are reconstructed using the limited slopes

$$u_i^{(m)}(x_{i+\frac{1}{2}}^-) = \bar{u}_i + \tilde{u}_i^m, \quad u_i^{(m)}(x_{i-\frac{1}{2}}^+) = \bar{u}_i - \tilde{\tilde{u}}_i^m. \quad (30)$$

For $k = 0, 1, 2$, the procedure uniquely reconstruct a polynomial of degree $k$. For $k \geq 3$ there is more freedom since the cell average and the two trace values do not completely determine the polynomial. Here we follow a simple approach given by [6]:

1. If the limiter does not modify the trace values, i.e.,

$$u_i^m(x_{i-\frac{1}{2}}^+) = u_i(x_{i-\frac{1}{2}}^+), \quad u_i^m(x_{i+\frac{1}{2}}^-) = u_i(x_{i+\frac{1}{2}}^-), \quad (31)$$

   then we take $u_i^m(x) = u_i(x)$.
2. Otherwise, let $u_i^1(x) \in P^1(E_i)$ be the $L^2$ projection of $u_i(x)|_{E_i}$. Take $u_i^{(m)}(x)|_{E_i} = \Pi_h(u_i^1(x))$, where $\Pi_h$ denotes the limiter.

For 2D problems, this limiting process can be carried out by applying the limiter sequentially in the $x$- and $y$- directions.

4.2 Zhong and Shu's simple WENO limiter

In this subsection we give a brief description of the simple WENO limiter developed by Zhong and Shu [40] as follows.

In one dimensional space, let us denote the DG solutions on the cells $E_i$, $E_{i-1}$, $E_{i+1}$ as $u_0(x), u_1(x), u_2(x)$, respectively. In order to ensure the reconstructed polynomial maintains the original cell average of $u_0(x)$ in the target cell $E_i$, the following modifications are made:

$$\tilde{u}_1(x) = u_1(x) - \bar{u}_1 + \bar{u}_0, \quad \tilde{u}_2(x) = u_2(x) - \bar{u}_2(x) + \bar{u}_1(x), \quad (32)$$

where

$$\bar{u}_0 = \frac{1}{|E_i|}\int_{E_i} u_0(x)\mathrm{d}x, \quad \bar{u}_1 = \frac{1}{|E_i|}\int_{E_i} u_1(x)\mathrm{d}x, \quad \bar{u}_2 = \frac{1}{|E_i|}\int_{E_i} u_2(x)\mathrm{d}x. \quad (33)$$

The final nonlinear WENO reconstruction polynomial $u_0^{new}(x)$ is defined by a convex combination of these modified polynomials:

$$u_0^{new}(x) = \omega_0 \tilde{u}_0(x) + \omega_1 \tilde{u}_1(x) + \omega_2 \tilde{u}_2(x). \quad (34)$$

The normalized nonlinear weights $\omega_j$ are defined as

$$\omega_j = \frac{\bar{\omega}_j}{\sum_l \bar{\omega}_l}, \quad (35)$$

where the non-normalized nonlinear weights $\bar{\omega}_j$ are functions of the linear weights $\gamma_j$ and the so-called smoothness indicators $\beta_j$:

$$\bar{\omega}_j = \frac{\gamma_j}{(\epsilon_0 + \beta_l)^r}. \quad (36)$$

In [40], the smoothness indicator for the 1D case is taken as:

$$\beta_j = \sum_{l=1}^{k} \int_{E_i} |E_i|^{2l-1} \big(\frac{\partial^l}{\partial x^l} u_i(x)\big)^2 \mathrm{d}x. \tag{37}$$

In the numerical simulations, we choose $\epsilon_0 = 10^{-6}$, $r = 2$ and the linear weights taken as

$$\gamma_0 = 0.998, \quad \gamma_1 = 0.001, \quad \gamma_2 = 0.001. \tag{38}$$

In [40], the readers could find more details about the simple WENO limiter for the two-dimensional space.

4.3 The high order shock-capturing moment limiter

The limiting strategy proposed by Moe et al. [23] can be viewed as a novel extension of the finite volume Barth-Jespersen limiter [3] and the modern maximum principle preserving DG schemes developed by Zhang and Shu[38]. We give a brief description of the implementation of this momentum limiter. Readers are referred to [23] for more details. The basic procedure consists of the following steps.

**Step 1.** For each mesh element $E_i$, compute an approximate maximum and minimum of the solution $u_i(x)$:

$$u_{M_i} := \max_{x \in \mathcal{X}_i}\{u_i(x)|_{E_i}\}, \quad u_{m_i} := \min_{x \in \mathcal{X}_i}\{u_i(x)|_{E_i}\}. \tag{39}$$

where the points set $\mathcal{X}_i$ consists of Gaussian quadrature points as well as corner points and quadrature points along the element boundaries.

**Step 2.** Consider the set $N_{E_i}$ of all neighbors of $E_i$ excluding $E_i$ itself, and compute an approximate upper and lower bounds:

$$M_i := \max\{\bar{u}_i + \alpha(h), \max_{j \in N_{E_i}}\{u_{M_j}\}\}, \tag{40}$$

$$m_i := \min\{\bar{u}_i - \alpha(h), \min_{j \in N_{E_i}}\{u_{m_j}\}\}, \tag{41}$$

The scalar function $\alpha(h) = \alpha h^{1.5} \geq 0$ is a tolerance function that depends on the problem. When $\alpha(h) = 0$, smooth extrema will be clipped, when $\alpha(h)$ goes to zero slow enough, the limiter will eventually turn off and not clip smooth extrema.

**Step 3.** Define

$$\theta_{M_i} := \Phi\big(\frac{M_i - \bar{u}_i}{u_{M_i} - \bar{u}_i + \epsilon_1}\big), \quad \theta_{m_i} := \Phi\big(\frac{m_i - \bar{u}_i}{u_{m_i} - \bar{u}_i - \epsilon_1}\big), \tag{42}$$

where $\epsilon_1 = 10^{-6}$ is a small number to prevent the denominator from becoming zero and the cut off function $0 \leq \Phi(s) \leq 1$ is

$$\Phi(s) := \min\{\frac{s}{1.1}, 1\}. \tag{43}$$

**Step 4.** Define the rescaling parameter as

$$\theta_i := \min\{1, \theta_{m_i}, \theta_{M_i}\}. \tag{44}$$

**Step 5.** Finally, rescale the approximate solution on the element $E_i$ as

$$u_i^{new}|_{E_i} := \bar{u}_i + \theta_i\big(u_i(x)|_{E_i} - \bar{u}_i\big). \tag{45}$$

4.4 Detector

The limiting procedures of the WENO limiter and Moe's limiter will be costly, in order to reduce the computation cost, we do not need to apply the limiters to every element, instead we first identify the troubled cells which might need the limiting procedure. Then we apply limiters only to those cells.

In our work, the minmod-TVB limiter with $M_{TVB} = 0$ will be used as a troubled cell indicator. When the limiter returns something else in the first argument, the cell is indicated as a troubled cell and will be limited.

**Table 3** Space accuracy tests of (46) at $T = 0.75/\pi$ ($\Delta t = 0.0005$).

| $k$ | $N$ | Without limiter | | TVB | | WENO | | Moe $\alpha = 0$ | | Moe $\alpha = 0.5$ | |
|---|---|---|---|---|---|---|---|---|---|---|---|
| | | $L^2$ error | order | $L^2$ error | order | $L^2$ error | order | $L^2$ error | order | $L^2$ error | order |
| 1 | 40  | 3.212e-03 | -- | 2.738e-02 | -- | 5.210e-02 | -- | 6.829e-03 | -- | 3.831e-03 | -- |
|   | 80  | 1.007e-03 | 1.67 | 8.178e-03 | 1.74 | 2.246e-03 | 4.54 | 1.742e-03 | 1.97 | 1.058e-03 | 1.86 |
|   | 160 | 2.858e-04 | 1.82 | 2.363e-03 | 1.79 | 2.850e-04 | 2.88 | 4.342e-04 | 2.00 | 2.891e-04 | 1.87 |
|   | 320 | 7.703e-05 | 1.89 | 7.016e-04 | 1.75 | 7.508e-05 | 1.92 | 1.044e-04 | 2.06 | 7.704e-05 | 1.91 |
| 2 | 40  | 1.996e-04 | -- | 2.799e-02 | -- | 5.872e-04 | -- | 5.204e-03 | -- | 3.703e-04 | -- |
|   | 80  | 2.518e-05 | 2.99 | 8.474e-03 | 1.72 | 2.318e-04 | 1.34 | 1.117e-03 | 2.22 | 3.460e-05 | 3.42 |
|   | 160 | 4.756e-06 | 2.40 | 2.515e-03 | 1.75 | 3.456e-05 | 2.75 | 2.383e-04 | 2.22 | 4.756e-06 | 2.86 |
|   | 320 | 7.822e-07 | 2.60 | 7.724e-04 | 1.70 | 3.706e-06 | 3.22 | 5.545e-05 | 2.11 | 7.822e-07 | 2.60 |
| 3 | 40  | 1.401e-05 | -- | 2.747e-02 | -- | 7.853e-04 | -- | 6.618e-03 | -- | 7.321e-04 | -- |
|   | 80  | 2.301e-06 | 2.61 | 8.158e-03 | 1.75 | 7.864e-05 | 3.32 | 1.347e-03 | 2.30 | 1.238e-05 | 5.89 |
|   | 160 | 2.576e-07 | 3.16 | 2.357e-03 | 1.79 | 3.390e-06 | 4.54 | 2.565e-04 | 2.38 | 2.576e-07 | 5.59 |
|   | 320 | 1.942e-08 | 3.73 | 7.210e-04 | 1.71 | 9.764e-08 | 5.12 | 5.534e-05 | 2.21 | 1.942e-08 | 3.73 |

**5 Numerical experiments**

In this section, we perform extensive numerical experiments to demonstrate the performance of the implicit-explicit Runge-Kutta discontinuous Galerkin method and the different limiters described in the previous sections. The program is implemented using the open source finite element library deal.ii [2]. First we present the numerical convergence orders, then we show the accuracy, efficiency and different features of the proposed scheme and limiters.

5.1 Accuracy tests

The accuracy tests are performed by solving the nonlinear Burgers equation

$$\begin{cases} \dfrac{\partial u}{\partial t} + \dfrac{\partial}{\partial x}(\dfrac{u^2}{2}) = 0, & x \in [0, 2], \\ u(0; t) = u(2; t), \end{cases} \tag{46}$$

with a continuous initial condition $u(x, t = 0) = \sin(\pi x)$. The exact solution to the Burgers equation satisfies $u - \sin(\pi(x - tu)) = 0$ and is smooth up to $t = \frac{1}{\pi}$.

For all accuracy tests, in the detecting step, we choose $M_{TVB} = 0$ for the minmod-TVB limiter. With this choice, many good cells are identified as troubled cells. Therefore, we can clearly see the effects of different limiters. When applying Moe's moment limiter, we choose two different values for $\alpha$, $\alpha = 0$ and $\alpha = 0.5$.

In Tables 3 and 4, we present the $L^2$-errors and convergence orders at $T = \frac{0.75}{\pi}$ and $\frac{1}{\pi}$ achieved by the DG scheme without limiter and with the TVB limiter, the WENO limiter and Moe's moment limiter. For smooth solutions at $T = \frac{0.75}{\pi}$, We can see that the TVB limiter can not keep the accuracy of the original scheme and reduces the order of high order schemes ($k = 2, 3$) to second order. The WENO limiter keeps both the designed order and the magnitude of errors of the original method. For Moe's moment limiter, the results obtained by different values of $\alpha$ are a bit different. When $\alpha = 0$, the moment limiter clips the smooth extrema (see the dashed cyan curve near the extrema in Fig. 2, therefore, reduces the original scheme to second order. When $\alpha$ is larger, the moment limiter does not clip the smooth extrema, thus keeps the designed order and the accuracy of the original scheme.

At $T = \frac{1}{\pi}$, the appearance of shock results in a large error near the discontinuity. Table 4 shows that all the schemes, with and without limiters, can not keep the high order accuracy and reduce the scheme to first order.

5.2 Numerical experiments in 1D

For the 1D simulations, we use three examples to illustrate the accuracy of the proposed DG scheme. Examples 1, 2 and 3 are taken from Ref. [19] and have been further studied by [37] using a moving mesh method.

In Eq. (9), we consider $F(u) = \frac{u^2}{u^2 + 0.5(1-u)^2}$, $\epsilon = 10^{-3}$ and initial condition

$$u(x; 0) = \begin{cases} 0, & x \in [0, 0.75], \\ u_0, & x \in (0.75, 2.25), \\ 0, & x \in [2.25, 3]. \end{cases} \tag{47}$$

**Table 4** Space accuracy tests of (46) at $T = 1.0/\pi$ ($\Delta t = 0.0005$).

| $k$ | $N$ | Without limiter | | TVB | | WENO | | Moe $\alpha = 0$ | | Moe $\alpha = 0.5$ | |
|---|---|---|---|---|---|---|---|---|---|---|---|
| | | $L^2$ error | order | $L^2$ error | order | $L^2$ error | order | $L^2$ error | order | $L^2$ error | order |
| 1 | 40 | 3.372e-02 | —— | 4.452e-02 | —— | 4.889e-02 | —— | 2.692e-02 | —— | 2.758e-02 | —— |
| | 80 | 1.867e-02 | 0.85 | 2.155e-02 | 1.05 | 1.685e-02 | 1.54 | 1.762e-02 | 0.61 | 1.752e-02 | 0.65 |
| | 160 | 1.039e-02 | 0.85 | 1.150e-02 | 0.91 | 8.086e-03 | 1.06 | 1.040e-02 | 0.76 | 1.039e-02 | 0.75 |
| | 320 | 5.802e-03 | 0.84 | 6.242e-03 | 0.88 | 4.298e-03 | 0.91 | 5.803e-03 | 0.84 | 5.802e-03 | 0.84 |
| 2 | 40 | 1.817e-02 | —— | 4.569e-02 | —— | 2.395e-02 | —— | 1.911e-02 | —— | 1.817e-02 | —— |
| | 80 | 1.011e-02 | 0.85 | 2.215e-02 | 1.04 | 1.009e-02 | 1.25 | 1.018e-02 | 0.91 | 1.011e-02 | 0.85 |
| | 160 | 5.642e-03 | 0.84 | 1.181e-02 | 0.91 | 4.708e-03 | 1.10 | 5.648e-03 | 0.85 | 5.642e-03 | 0.84 |
| | 320 | 3.1525-03 | 0.84 | 6.382e-03 | 0.89 | 2.281e-03 | 1.06 | 3.156e-03 | 0.84 | 3.155e-03 | 0.84 |
| 3 | 40 | 1.182e-02 | —— | 4.522e-02 | —— | 2.877e-02 | —— | 1.140e-02 | —— | 1.184e-02 | —— |
| | 80 | 6.586e-03 | 0.84 | 2.195e-02 | 1.04 | 1.500e-02 | 0.94 | 6.763e-03 | 1.05 | 6.586e-03 | 0.85 |
| | 160 | 3.680e-03 | 0.84 | 1.160e-02 | 0.92 | 8.152e-03 | 0.88 | 3.692e-03 | 0.87 | 3.681e-03 | 0.84 |
| | 320 | 2.060e-03 | 0.84 | 6.188e-03 | 0.91 | 4.481e-03 | 0.86 | 2.061e-03 | 0.84 | 2.060e-03 | 0.84 |

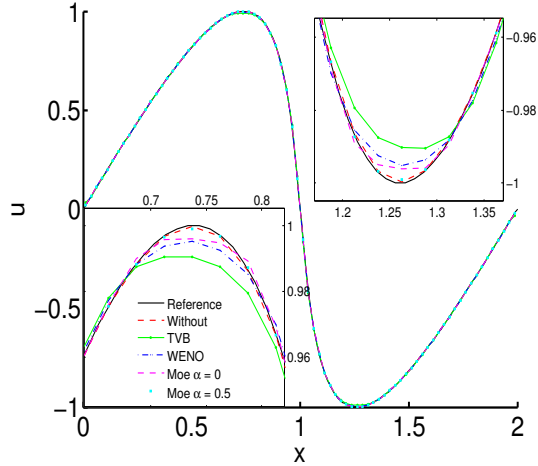

**Fig. 2** Solutions (one point per cell) computed using the DG scheme ($k = 1$) without limiter and with the TVB limiter, the WENO limiter, Moe's moment limiter ($\alpha = 0, 0.5$) for Burgers equation at $T = 0.75/\pi$.

With different combinations of $(\tau, u_0)$, the MBLE admits several types of solutions: rarefaction wave, admissible Lax shock and undercompressive shock. For Examples 1-3, we solve the 1D MBLE to $T = 0.5$.

**Example 1.** $(\tau, u_0) = (5.0, 0.66)$.

The TW results in [19, 37] show that, for the pair $(\tau, u_0) = (5.0, 0.66)$, the solution consists of a monotone basin of height $\underline{u} = 0.2027$ in the left part and a non-monotone plateau of height $\bar{u} = 0.7130$ in the right part.

In Fig. 3 we present the numerical solutions obtained by the DG scheme with and without limiters. We use $N = 501$ grid points in the space discretization and $k = 1$ for the basis polynomials. We check the accuracy of the obtained solutions by verifying the computed overshoot values. Fig. 3 (left) shows all three limiters suppress the overshoot in both the left basin region and the right plateau region. With respect to overshoot values, Moe's moment limiters ($\alpha = 0, 100, 200$) perform better than the WENO limiter and the minmod-TVB limiter. With the increase of $\alpha$ from 0 to 200, the limiting effects become weaker and weaker and the limited solutions become more and more close to the original solution.

Fig. 3 (right) shows the numerical convergence of Moe's moment limiter with $\alpha = 100$ for the basis polynomials with degree $k = 1, 2, 3$. With the increase of $k$, the solutions converge to the TW solution. The overshoot values obtained by the highest order DG scheme ($k = 3$) almost coincide with the predicated TW values. In Fig. 4, we plot the solutions computed using the forth order DG method with different choices of limiter. Again, Moe's moment limiter performs best, followed by the WENO limiter and the minmod-TVB limiter.

In the above simulations, for Moe's moment limiter with $\alpha = 100$ we only needs 501 grid points to accurately approximate the overshoot values. These solutions are comparable (in the eyeball norm) to the solutions in Fig. 4 in Ref. [19] obtained by the splitting method using 16384 points with the minmod-based reconstruction and 4096 points with the WENO5 reconstruction. From the above comparisons, we can conclude that in the simulation of MBLE, the DG method with Moe's moment limiter is superior to the splitting methods with minmod and WENO5 reconstruction.

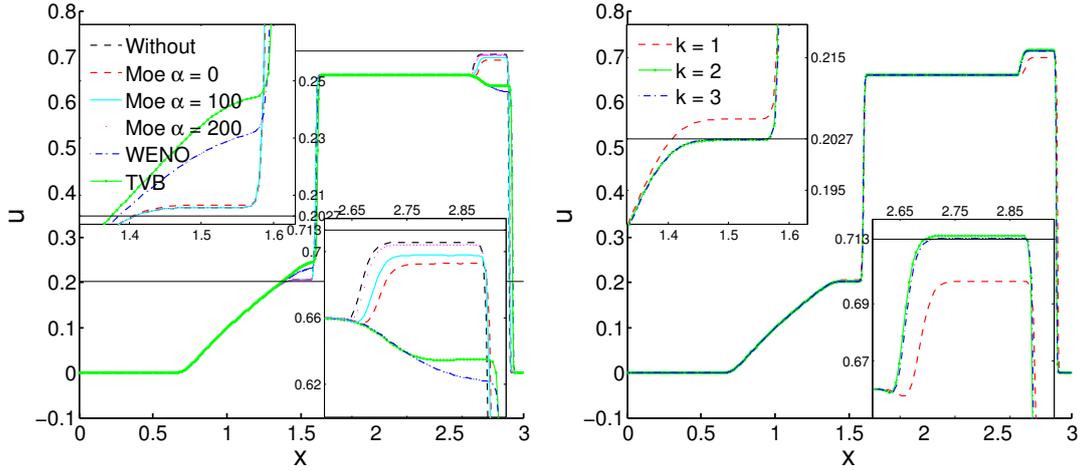

**Fig. 3** Example 1. Left: comparison between the DG solutions without limiter and with the TVB limiter $M_{TVB} = 0$, the WENO limiter and Moe's moment limiter ($\alpha = 0, 100, 200$) (in the upper left zoom-in figure, the curves of 'Without', 'Moe $\alpha = 100$', 'Moe $\alpha = 200$' coincide); right: comparison between the solutions using Moe's moment limiter ($\alpha = 100$) for $k = 1, 2, 3$.

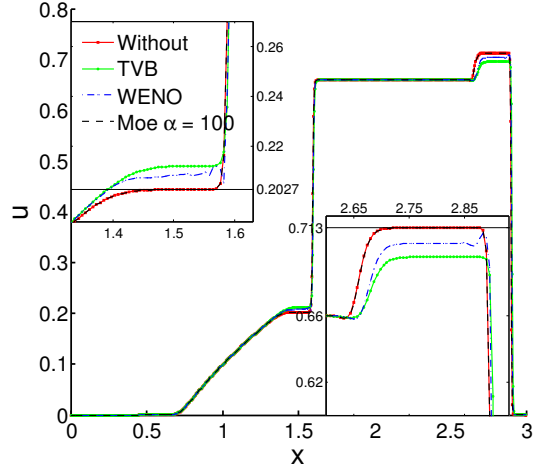

**Fig. 4** Example 1: solutions computed using the DG scheme without limiter and with the TVB limiter $M_{TVB} = 0$, the WENO limiter and Moe's moment limiter ($\alpha = 100$) for $k = 3$.

**Example 2.** $(\tau, u_0) = (5, 0.52)$.

In the second example, we decrease $u_0$ to 0.52. The bifurcation diagram shows this pair still admits a plateau of height $\bar{u} = 0.713$. Due to the decrease of $u_0$, the speed of the shock to the left of the plateau

$$s_2 = \frac{F(0.52) - F(0.713)}{0.52 - 0.713} \approx 1.1597, \tag{48}$$

is bigger than the shock speed in Example 1

$$s_1 = \frac{F(0.66) - F(0.713)}{0.66 - 0.713} \approx 0.7963, \tag{49}$$

Since the speed of the right undercompressive shock remains the same, we can expect a narrower plateau in the moving front.

We present the results computed by DG method with Moe's moment limiter in Fig. 5 (left). With the increase of polynomial degree, the overshoot plateau value becomes more and more close to the TW value.

**Example 3.** $(\tau, u_0) = (3.5, 0.85)$.

The third pair of $(\tau, u_0)$ corresponds to region $A_2$ in the bifurcation diagram in Fig. 1. The TW solution is different from those in Examples 1 and 2: the right part consists of a rarefaction wave connecting a plateau of height $\bar{u} = 0.6938$ trailing an undercompressive shock, the left part also consists of a rarefaction wave connection a basin of height $\underline{u} = 0.1036$ trailing an undercompressive shock. The zoomed-in figures in Fig. 5 (right) again illustrate the accuracy of high order method.

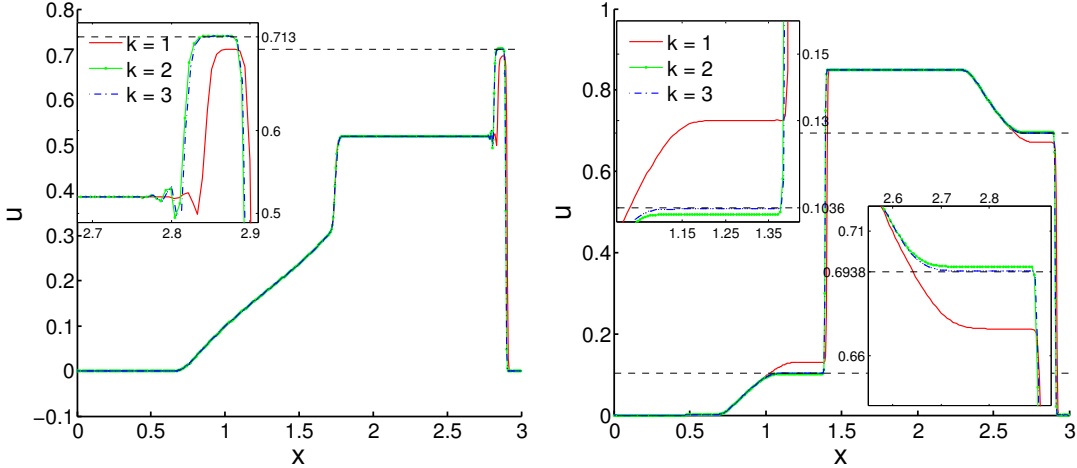

**Fig. 5** Example 2 (left) and Example 3 (right): solutions compute using the DG scheme with Moe's moment limiter ($\alpha = 100$) for $k = 1, 2, 3$.

5.3 Numerical experiments in 2D

In this subsection, we will study the effects of limiters in 2D simulations of the MBLE. In Ref. [37], the authors studied the 2D simulations of the MBLE by using a moving mesh finite difference method. Their results show that to compute accurate saturation plateaus, adaptive mesh gives a great advantage. In this work, we will show the advantage of using high order schemes.

**Example 4.** First, we consider the MBLE with $\tau = 0$. The functions and parameter are taken as

$$\begin{cases} F(u) = \dfrac{u^2}{u^2 + (1-u)^2}, \\ G(u) = F(u)(1 - 5(1-u)^2), \\ \epsilon = 0.01. \end{cases} \tag{50}$$

The initial data is

$$u(x, y, 0) = \begin{cases} 1, & x^2 + y^2 < 0.5, \\ 0, & \text{otherwise.} \end{cases}$$

We solve the equation in the square domain $[-1.5, 1.5] \times [-1.5, 1.5]$ to $T_{end} = 0.5$.

This test case is first proposed and solved in [20] and has been chosen as a benchmark test by many researchers [21, 39, 22, 36].

In Fig. 6 we present the numerical solutions and contour plots obtained by the DG method with and without limiters. It can be observed that on a mesh with $101^2$ points, oscillations will appear along the moving front when limiters are not applied to the DG solution. Different limiters also affect the solutions differently. The TVB limiter with $M_{TVB} = 50$ gives the most smooth solution and is comparable to those obtained in previous works [20, 39, 22, 36, 12]. When the WENO limiter or Moe's moment limiter ($\alpha = 10$) are used, the strong oscillations along the moving front are weakened, but slight fluctuations are still observable. We show the slice views at $y = 0.75$ of the solutions in Fig. 7 (left). We can clearly see the differences between the solutions. Without limiters, the DG scheme presents strong oscillations near the left top and right bottom corners, when limiters are activated, oscillations will be suppressed. The minmod-TVB limiter with $M_{TVB} = 50$ gives the most smooth solution, followed by the WENO limiter and Moe's moment limiter $\alpha = 10$.

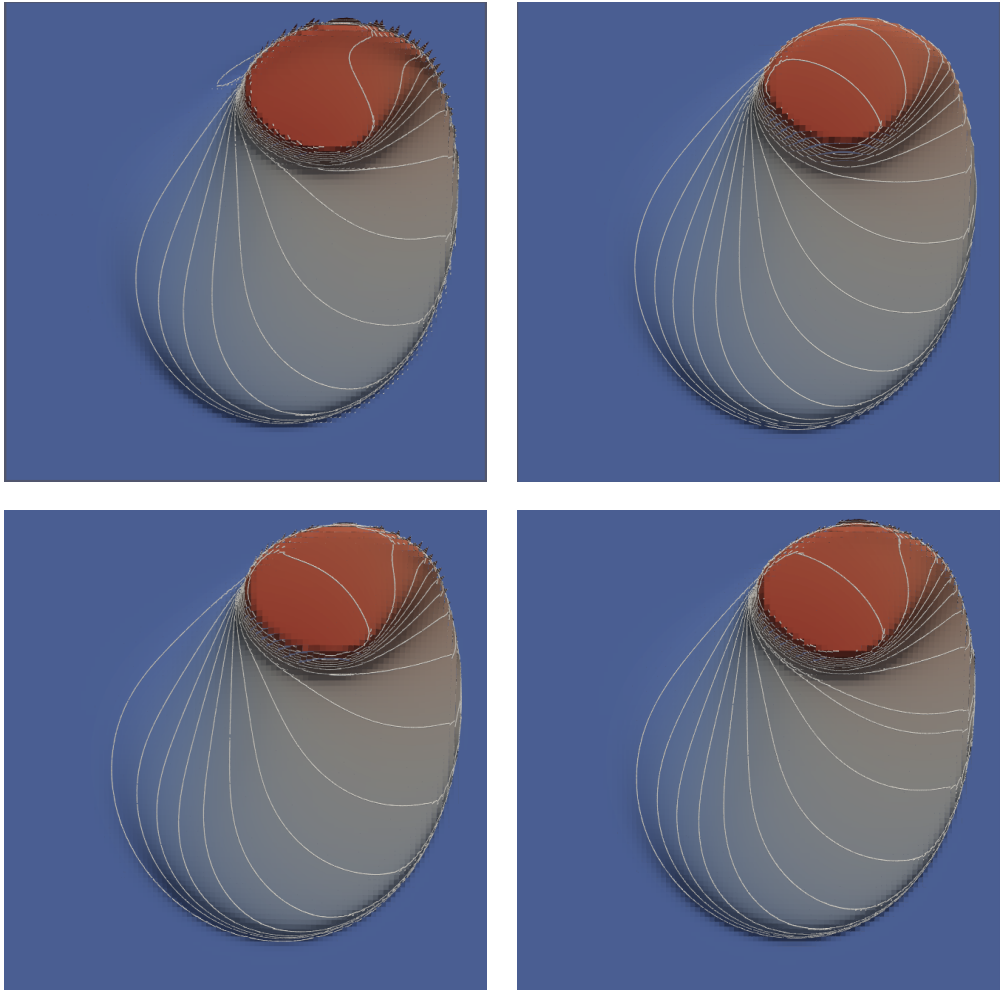

**Fig. 6** Example 4: solutions and contour plots (20 contour lines in $[0, 1]$) computed using the DG scheme without limiter (top left) and with the TVB limiter $M_{TVB} = 50$ (top right), the WENO limiter (bottom left) and Moe's moment limiter $\alpha = 10$ (bottom right).

**Example 5.** We study two different initial conditions of the 2D MBLE with $\tau = 0.5$. The functions and parameters are the same as those used in Example 4. The initial conditions are taken as: one with a cylindrical shape

$$\text{Example 5-1:} \quad u(x, y, 0) = \begin{cases} 0.9, & x^2 + y^2 < 0.5, \\ 0, & \text{otherwise,} \end{cases} \quad (x, y) \in [-1.5, 1.5] \times [-1.5, 1.5], \quad (51)$$

and one with a cubic shape

$$\text{Example 5-2:} \quad u(x, y, 0) = \begin{cases} 0.9, & x^2 < 0.5, y^2 < 0.5, \\ 0, & \text{otherwise,} \end{cases} \quad (x, y) \in [-1.5, 1.5] \times [-1.5, 1.5]. \quad (52)$$

When dynamic coefficient $\tau$ is not zero, we use the TW analysis in Section 2.2 to predict the behavior of the solution. The results in [36] shows that with $\tau = 0.5$, a saturation plateau of height $\bar{u} \approx 0.97$ will be developed near the shock front in the $y$-direction. Fig. 8 gives the 3D view of the numerical solutions computed at $T = 0.5$ by the DG scheme without limiter and with the TVB limiter, the WENO limiter and

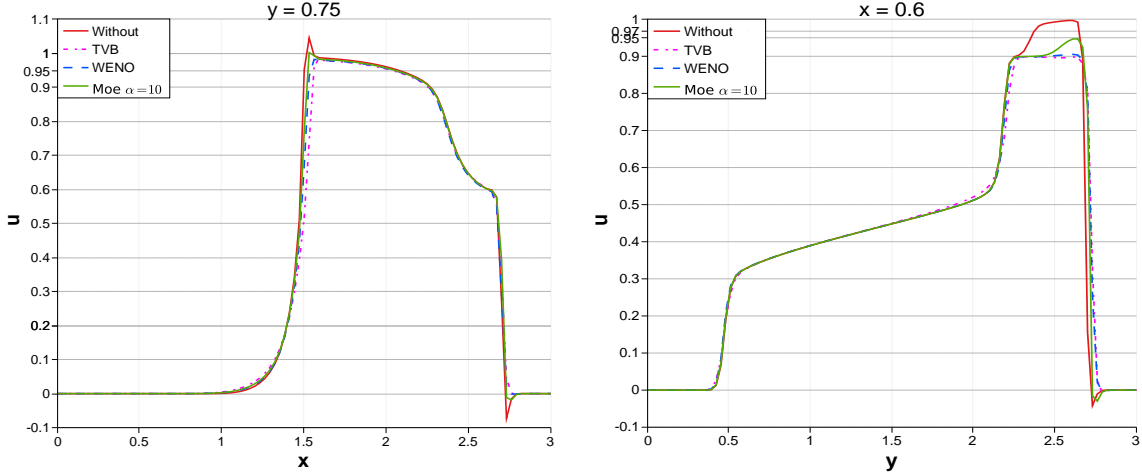

**Fig. 7** Left: slices at $y = 0.75$ of DG solutions computed without limiter, with the TVB limiter $M_{TVB} = 50$, the WENO limiter and Moe's moment limiter $\alpha = 10$ for Example 4. Right: slices at $x = 0.6$ of DG solutions computed without limiter, with the TVB limiter $M_{TVB} = 50$, the WENO limiter and Moe's moment limiter $\alpha = 10$ for Example 5-1.

Moe's moment limiter on a mesh with $101^2$ points. The TVB limiter with $M_{TVB} = 50$ almost smooths out the plateau while the WENO limiter gives a slight higher plateau and Moe's moment limiter obtained the highest plateau value as also shown in Fig. 7 (right). From Figs. 7 and 8, we can conclude that Moe's moment limiter with $\alpha = 10$ is the least diffusive among all three limiters.

Since Moe's moment limiter performs best in this example, next we will solve this problem with different polynomial orders. Fig. 9 shows the 3D views and slice views for $k = 1, 2, 3$. When the polynomial degree increases from $k = 1$ to $k = 3$, the overshoot plateau value become more and more close to the TW value.

Notice that in Ref. [36], the authors carried out the simulations on a uniform mesh with $1000^2$ points and an adaptive mesh with $300^2$ grid points. Our results on a mesh with $101^2$ points again show the advantages of using high order DG scheme and suitable limiters.

In the final test, we show the results computed using the cubic shape-initial condition on a mesh with $201^2$ grid points in Fig. 10. Similar to the previous results in Ref. [36], the non-monotone plateaus are located near the shock front in the $y$-direction and become thinner and lower along the positive $x$-direction because of the rarefaction waves created by the flux in the $x$-direction. Again, the high order DG scheme shows its advantage over lower scheme.

## 6 Conclusions

In this work, we implement a symmetric interior penalty discontinuous Galerkin scheme for the modified Buckley-Leverett equation in both 1D and 2D spaces. The minmod-TVB limiter, the WENO limiter and Moe's moment limiter were implemented to limiter the DG solutions. Specifically, comparing with the TVB limiter and the WENO limiter, the DG method with Moe's moment limiter yielded very little numerical diffusion and outperformed the other two limiters in approximating the TW solutions. In addition, our results also demonstrated the advantage of using high order DG method: on a 1D mesh with 501 grid points, the DG method with Moe's moment limiter easily outperformed the operator splitting method with second order minmod limiter (16384 points) and the fifth order WENO5 limiter (4096 points) in Ref. [19].

**Acknowledgements**

H. Zhang gratefully acknowledges the financial support from the China Scholarship Council (No. 201503170430). Y. Guo is supported by the National Natural Science Foundation of China (No. 11501570, 91530106). W. Li is supported by the National Key Research and Development Program of China (No. 2017YFF0210701).

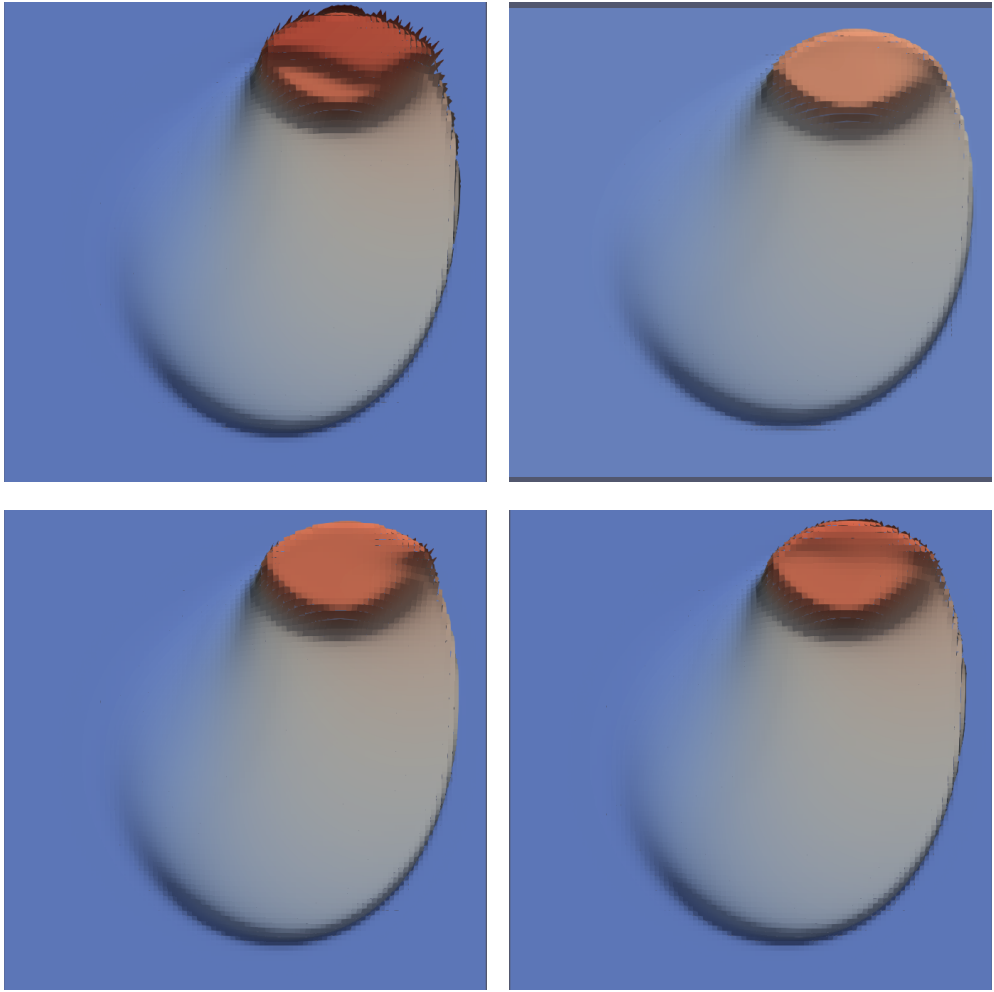

**Fig. 8** Example 5-1. Solutions computed using the DG scheme ($k = 1$) without limiter (top left) and with the TVB limiter $M_{TVB} = 50$ (top right), the WENO limiter (bottom left) and Moe's moment limiter $\alpha = 10$ (bottom right).

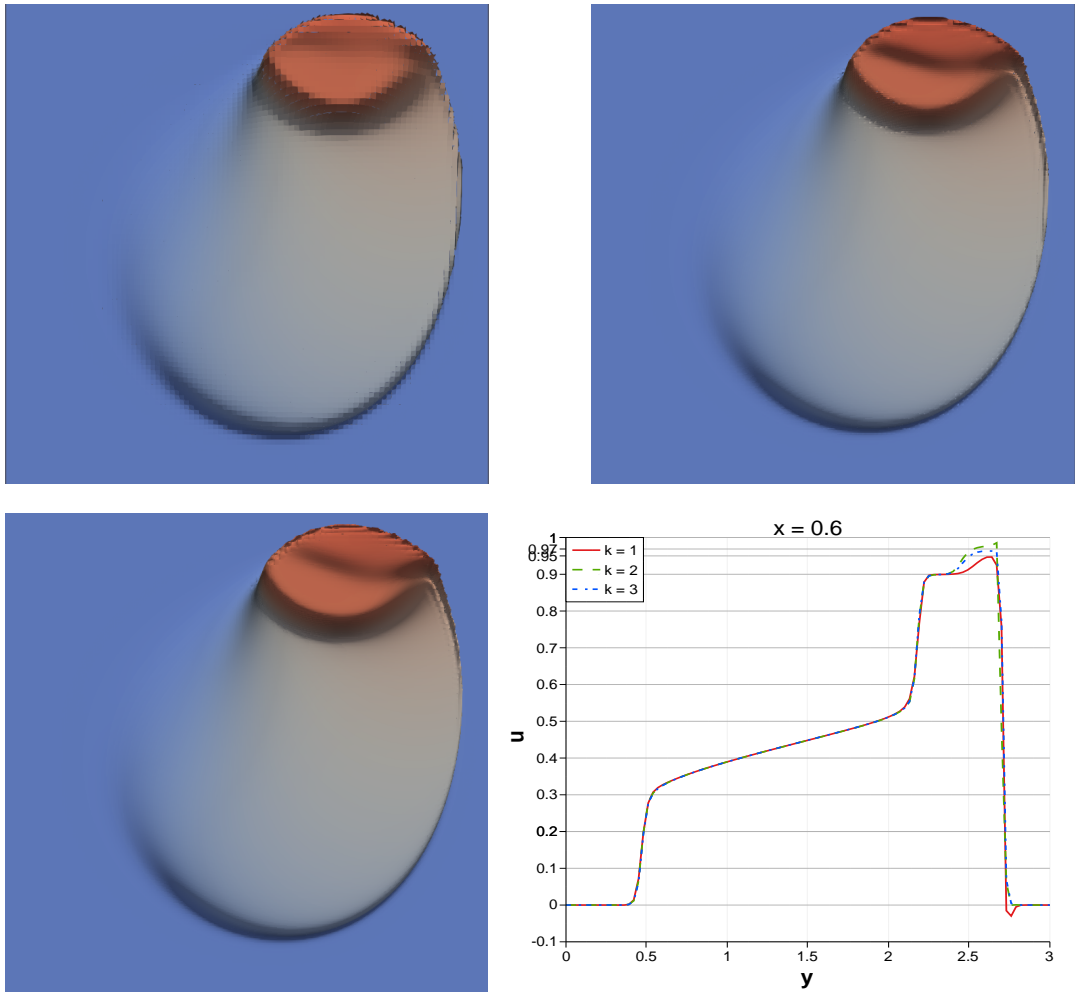

**Fig. 9** Example 5-1. Solutions computed using the DG scheme with Moe's moment limiter $\alpha = 10$ for $k = 1$ (top left), $k = 2$ (top right) and $k = 3$ (bottom left). Slices at $x = 0.6$ for $k = 1, 2, 3$ (bottom right).

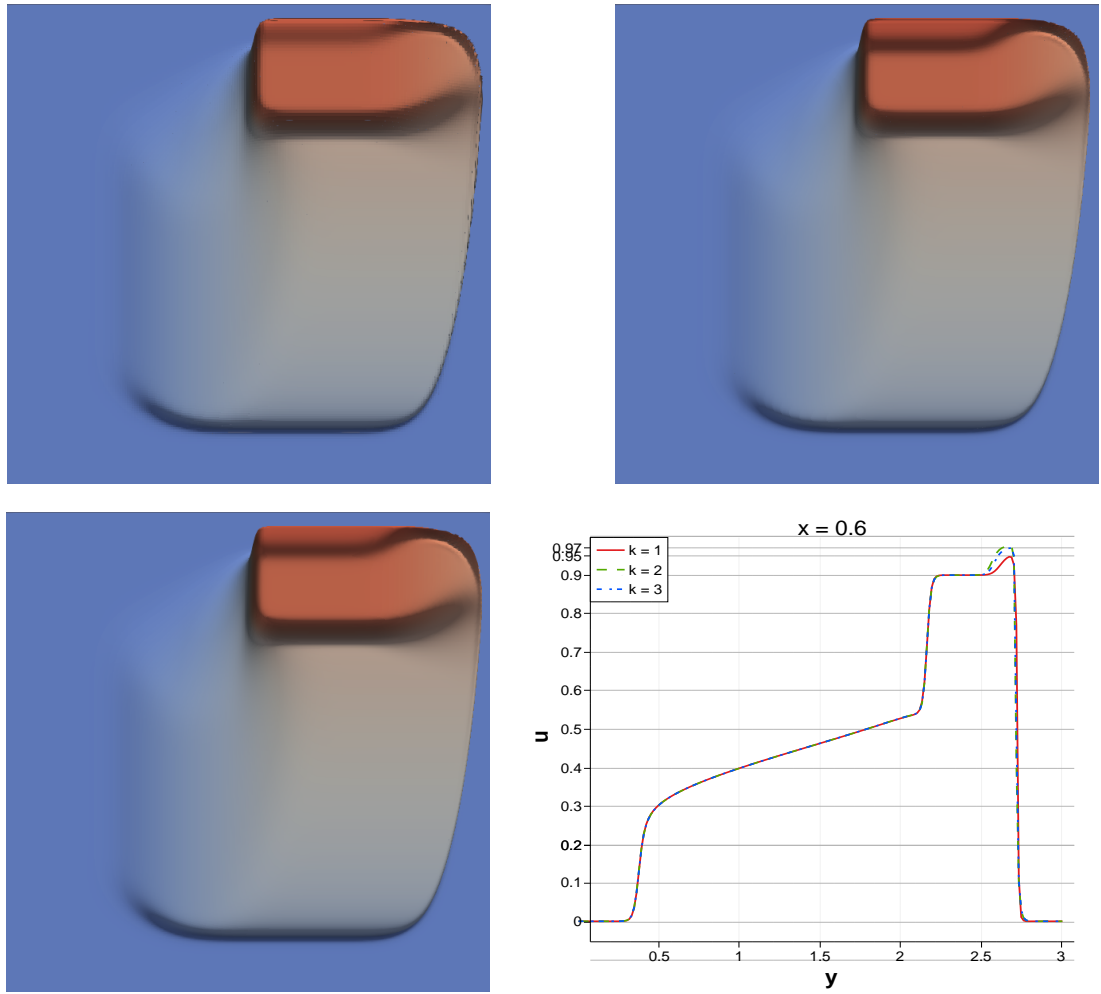

**Fig. 10** Example 5-2. Solutions computed using the DG scheme with Moe's moment limiter $\alpha = 10$ for $k = 1$ (top left), $k = 2$ (top right) and $k = 3$ (bottom left). Slices at $x = 0.6$ for $k = 1, 2, 3$ (bottom right).